\documentclass[12pt]{article}

\newcommand{\nc}{\newcommand}
\nc{\ov}{\over}
\nc{\iy}{\infty}
\nc{\s}{\sigma}
\nc{\al}{\alpha}
\nc{\x}{\xi}
\nc{\inv}{^{-1}}
\renewcommand{\sp}{\vspace{1ex}}
\nc{\be}{\begin{equation}} \nc{\ee}{\end{equation}}
\nc{\bqn}{\begin{eqnarray}} \nc{\eqn}{\end{eqnarray}}
\nc{\nn}{\nonumber}
\nc{\ph}{\varphi} 
\renewcommand{\th}{\theta}
\nc{\ch}{\raisebox{.4ex}{$\chi$}} 
\nc{\dl}{\delta}
\nc{\ep}{\varepsilon} 
\nc{\e}{\eta} 
\nc{\tl}{\widetilde}
\nc{\twotwo}[4]{\left(\begin{array}{cc}#1&#2\\&\\#3&#4\end{array}\right)}
\nc{\Z}{\mathbf{Z}}
\nc{\Zp}{{\Z^+}}
\nc{\bR}{\mathbf{R}}
\nc{\C}{\mathbf{C}}
\nc{\T}{\mathbf{T}}

\begin{document}

\begin{center}{\bf On the Determinant of a Certain Wiener-Hopf + Hankel Operator}\end{center}
\begin{center}{{\bf Estelle L. Basor}\\
{\it Department of Mathematics\\
California Polytechnic State University\\
San Luis Obispo, CA 93407, USA}}\end{center}

\begin{center}{{\bf Torsten Ehrhardt}\\
{\it Fakult\"at f\"ur Mathematik\\
Technische Universit\"at Chemnitz\\
09107 Chemnitz, Germany}}\end{center}

\begin{center}{{\bf Harold Widom}\\
{\it Department of Mathematics\\
University of California\\
Santa Cruz, CA 95064, USA}}\end{center}

\sp

\begin{abstract} We establish an asymptotic formula for determinants of truncated 
Wiener-Hopf\linebreak+Hankel operators with symbol equal to the exponential of a constant
times the characteristic function of an interval. This is done by reducing it to 
the corresponding (known) asymptotics for truncated Toeplitz+Hankel operators. 
The determinants in question
arise in random matrix theory in determining the limiting distribution
for the number of eigenvalues in an interval for a scaled Laguerre ensemble
of positive Hermitian matrices. 
\end{abstract}

\begin{center}{\bf I. Introduction}\end{center}\sp

It has long been known that the asymptotic evaluation of Fredholm determinants of truncated
Wiener-Hopf operators $W_R(\s)$ is important in the study of certain random matrix problems. 
These operators act on $L^2(0,R)$ according to the rule
$$
f(x)\mapsto g(x)=f(x)+\int_{0}^{R}k(x-y)f(y)\,dy
$$
with the kernel of the integral operator given by
$$k(x) = {1\ov 2\pi}\int_{-\iy}^{\iy} ( \s(\xi) -1)\,e^{-ix\xi}\,d\xi.$$
The symbol $\s(\xi)$ is defined on the real line $\bR$ and it is assumed that 
$\s(\xi)-1\in L^1(\bR)$.  If we think of
$\lambda_{1}, \dots , \lambda_{N}$ as the eigenvalues of a random Hermitian matrix in the
classical Gaussian Unitary Ensemble
(GUE) the Fourier transform of the distribution function for the random variable given by
$$\sum_{i=1}^{N}f(\lambda_{i}/\sqrt{2N})$$ can be shown to be equal (after an appropriate
scaling limit) to the Fredholm determinant
$$\det W_{R}(\s)$$
 with symbol
$$\s(\xi)=e^{i\alpha f(\x)}.$$ An application of a Szeg\"{o} type theorem, which computes these
determinants asymptotically, shows that the random variable is asymptotically normal. This holds
if the function $f$ satisfies the conditions required by
the Szeg\"{o}  theorem. In particular, the symbol should not have jump discontinuities.
In the case where the function $f$ is a characteristic function the symbol has
at least two jumps and the classical results do not apply. However, generalizations of the theorem
starting with \cite{BW2} and mostly recently with \cite{BW} have been obtained so that one can
now compute the asymptotics of these determinants for symbols with jumps and also
other singularities of a  certain type. For a history of this problem, the reader is referred to
\cite{Arch,BS}. It is interesting to note here that the asymptotics for the case where jumps occur have a different form that make
the associated distributions not asymptotically normal. For more information about the connection between the random variables and Fredholm determinants we refer to \cite{Ba}.

In the case of other random matrix ensembles other operators arise. In particular, for 
the so-called Laguerre ensembles of positive random Hermitian matrices, the operators of 
interest are not Wiener-Hopf, but Bessel operators, where the Fourier transform 
is replaced by the Hankel transform. In this case a Szeg\"{o} type theorem was proved 
as well for symbols that are smooth. These results are contained in \cite{BE1}. 

In two important cases, the Bessel operators reduce to operators which are a sum or
difference of truncated Wiener-Hopf  and Hankel operators. These operators also act on
$L^{2}(0, R)$ and according to the rule
$$
f(x)\mapsto g(x)=f(x)+\int_{0}^R\Big(k(x-y)\pm k(x+y)\Big)f(y)\,dy,
$$
where $k$ is the same as above. We will denote them by 
$$W_{R}(\s) \pm H_{R}(\s).$$
Generalizations of Szeg\"o type theorems for symbols that are not smooth have
been very difficult to obtain for this class of operators. One of the main difficulties is that
while one can attempt to apply localization techniques which reduce the problem of singularities
to a special family of symbols, the determinants for this family are not explicitly known.

It was exactly the analogous problem which delayed the generalization of the Szeg\"{o} theorem
in the Wiener-Hopf case to singular symbols. However, the difficulty was overcome in \cite{BW}
by finding alternative Fredholm determinant representations for the determinants of the truncated
Wiener-Hopf operators and then comparing them to
corresponding formulas for determinants of finite Toeplitz matrices. These matrices are defined
by
\[T_n(\ph)=(\ph_{j-k})_{j,k=0,\dots,n-1}.\]
Here $\ph_k$ stands for the $k$-th  Fourier coefficients of a function $\ph$ defined on the
unit circle. 
For these the asymptotics were, conveniently, known. In \cite{BW} it was shown that
for certain classes of singular symbols the
Toeplitz and Wiener-Hopf determinants were (up to a simple constant factor) asymptotically equal
when $R \sim 2n$.
The alternative representations were obtained by using an identity of Borodin and Okounkov
for Toeplitz determinants \cite{BO} and its Wiener-Hopf analogue \cite{BC}. These formulas were only valid
for smooth symbols and so it was necessary to introduce a parameter to
regularize the symbols, apply the identity
and then take a limit. This gave exact formulas for the determinants for the
singular symbols. The final step was to prove that these were asymptotically equal.

In this paper we follow a similar path to find analogous asymptotic formulas for the
truncated Wiener-Hopf + Hankel operators
\[W_{R}(\s) + H_{R}(\s)\]
in a case where $\s$ is the characteristic function of an interval and so has two jump
discontinuities. Fortunately, the asymptotics for the corresponding Toeplitz + Hankel matrices
\[T_n(\ph)+H_n(\ph)=((\ph_{j-k}) + (\ph_{j+k+1}))_{j,k=0,\dots,n-1}\]
are known for these symbols \cite{BE2} and thus the comparison can be made.

The symbol of interest is
\[\s(\x)=\left\{\begin{array}{ll} 1&{\rm if}\ |\x|>1,\\&\\e^{-2\pi i\al}&{\rm if}\ |\x|<1,
\end{array}\right.\]
and the corresponding kernel is
\[(e^{-2\pi i\al}-1)\left({\sin\,(x-y)\ov \pi\,(x-y)}+{\sin\,(x+y)\ov \pi\,(x+y)}\right).\]
\sp

\noindent{\bf Theorem}. With $\s$ as above we have, for $|\Re\,\al|<1/2$,
\[\det\,(W_{R}(\s) + H_{R}(\s))\sim e^{-2i\al R}\,R^{-3\al^2}\,2^{4\al^2}\,G(1-2\al)\,G(1+2\al),\]
where $G$ is the Barnes $G$-function. \sp

Here is an outline of the paper. In the next section we derive Fredholm determinant representations for the determinants of
$T_n(\ph)+H_n(\ph)$ and $W_{R}(\s) + H_{R}(\s)$ in case of smooth symbols. In Section~III we
introduce a parameter to regularize our symbols and apply the formulas. We then take quotients and
let the parameter tend to
its limit, obtaining in this way an exact formula for the quotient of the determinants with the
singular symbols. In Section~IV we show that (up to a simple factor) they are asymptotically equal.

The stated results yields asymptotic information about the random variable
$$\sum_{i=1}^{N}f(\lambda_{i}/\sqrt{2N})$$ where $f$ is the characteristic function of the interval
$(0,1).$ This random variable counts the number of eigenvalues in the interval
$(0, \sqrt{2N})$. The techniques here can likely be modified to apply to other symbols as
well so that more general discontinuous symbols can be also studied.

\begin{center}{\bf II. Formulas for truncated determinants}\end{center}\sp

If $A$ is an invertible operator on Hilbert space of the form identity$+$trace class then for
projections $P$ and $Q=I-P$ we have
\be\det\,PAP=(\det\,A)\cdot(\det\,QA\inv Q).\label{detdet}\ee
This can be used to derive a formula of Borodin and  Okounkov \cite{BO} for the
finite Toeplitz determinants $\det\,T_n(\ph)$  and an anologous
formula for the finite Wiener-Hopf opertors $W_R(\s)$ \cite{BC}. For the former $P$ is
the projection $P_n$ onto the subspace of $\ell^2(\Zp)$, $\Zp=\{0,1,\dots,\}$,
spanned by $z^k$ ($k=0,\dots,n-1$), and for the latter $P$ is the 
projection onto the subspace $L^2(0,R)$ of $L^2(0,\iy)$.
A proof of the Borodin Okounkov formula for Toeplitz determinants based on 
formula (\ref{detdet}) implicitly is in \cite{BW} and explicitly is in \cite{Bo}.

Here we shall obtain analogous formulas for the determinants of 
the operators $T_n(\ph)+H_n(\ph)$ and $W_R(\s)+H_R(\s)$. We 
assume about $\ph$ and $\s$ that they have sufficiently smooth 
logarithms (with $\log\s$ also belonging to $L^1(\bR)$). Thus 
they have factorizations $\ph(z)=\ph^+(z)\,\ph^-(z)$ and 
$\s(\x)=\s^+(\x)\,\s^-(\x)$, with factors extending analytically 
in the usual way to be analytic and nonzero inside and outside the 
unit circle and above and below the real line, respectively. Moreover, we 
will restrict ourselves to the case where the symbols are even, i.e., 
$\ph(z\inv)=\ph(z)$ and $\s(-\x)=\s(\x)$, since that is all we shall need later. 
In this case we can assume without loss of generality that the factors 
of the factorizations  are related by $\ph^-(z)=\ph^+(z\inv)$ and 
$\s^-(\xi)=\s^+(-\xi)$. Formulas can also be obtained for non-even symbols, 
but they are more complicated.

In order to establish these formulas, let us recall that the Wiener-Hopf and 
the Hankel operators acting on $L^2(0,\iy)$ are defined by
$$
W(\sigma):f(x)\mapsto g(x)=f(x)+\int_{0}^\iy k(x-y) f(y)\,dy,
$$
$$
H(\sigma):f(x)\mapsto g(x)=\int_{0}^\iy k(x+y) f(y)\,dy,
$$
where $k(x)$ is given as before and $\s-1\in L^1(\bR)$. 
The Toeplitz and the Hankel operator acting on $\ell^2(\Zp)$ are given in terms
of their matrix representations  by
$$
T(\ph)=(\ph_{j-k})_{j,k\in\Zp},\qquad
H(\ph)=(\ph_{j+k+1})_{j,k\in\Zp},
$$
where $\ph_k$ are the Fourier coefficients of the function $\ph$ defined on 
the unit circle.

If $\ph$ and $\s$ satisfy the above mentioned assumption and are even, then the 
identities are
\be\det\,(T_n(\ph)+H_n(\ph))=G[\ph]^n\,\det\,[T(\ph\inv)\,(T(\ph)+H(\ph))]\,
\det\,(I+Q_n\,H(\ph_-/\ph_+)\,Q_n),\label{T+H}\ee
\be\det\,(W_R(\s)+H_R(\s))=G[\s]^R\,\det\,[W(\s\inv)\,(W(\s)+H(\s))]\,\det\,
(I+Q_R\,H(\s_-/\s_+)\,Q_R).\label{W+H}\ee
The constants $G[\ph]$ and $G[\s]$ are well known constants that 
appear also in the Szeg\"o type theorems,
$$
G[\ph]= \exp\left(\frac{1}{2\pi}\int_{0}^{2\pi}\log\ph(e^{ix})\,dx\right),
$$
$$
G[\s]= \exp\left(\frac{1}{2\pi}\int_{-\iy}^{\iy}\log\s(\xi)\,d\xi\right).
$$

To prove (\ref{T+H}) we write
\bqn
\det\,(T_n(\ph)+H_n(\ph))
&=&  \det\,P_n(T(\ph)+H(\ph))\,P_n \nn\\
&=&  \det\,P_n\,T(\ph^+)\,P_n\,T\inv({\ph^+})\,
(T(\ph)+H(\ph))\,T\inv({\ph^-})\,P_n\,T(\ph^-)\,P_n \nn\\
&=& G[\ph]^n\,\det\,P_n\,T\inv({\ph^+})\,
(T(\ph)+H(\ph))\,T\inv({\ph^-})\,P_n.\nn
\eqn
We apply (\ref{detdet}) with $P=P_n$ and $A=T\inv({\ph^+})\,(T(\ph)+H(\ph))\,T\inv({\ph^-})$.
Then $\det\,A$ equals the second factor on the right side of (\ref{T+H}) since
$T\inv({\ph^-})\,T\inv({\ph^+})=T(\ph\inv)$, and it remains to identify $A\inv$.
Of course
\[A\inv=T(\ph^-)\,(T(\ph)+H(\ph))\inv\,T(\ph^+).\]
The pleasant fact is that if $\ph$ is even and invertible then
\[(T(\ph)+H(\ph))\inv=T(\ph\inv)+H(\ph\inv).\]
(See \cite{BE}, Sec.~2. For non-even $\ph$ there is a more complicated formula.) Since
$T(\ph\inv)=T\inv({\ph^-})\,T\inv({\ph^+})$ this gives
\[A\inv=I+T(\ph^-)\,H(\ph\inv)\,T(\ph^+).\]
We now use twice the general identity
\[H(\psi_1\,\psi_2)=T(\psi_1)\,H(\psi_2)+H(\psi_1)\,T(\tl{\psi_2})\]
(where $\tl\psi(z)=\psi(z\inv))$ to obtain
\[T(\ph^-)\,H(\ph\inv)\,T(\ph^+)=T(\ph^-)\,H(\ph\inv\,\tl{\ph^+})=H(\ph^-\,\ph\inv\,\tl{\ph^+})=
H(\ph^-/{\ph^+})\]
since we have assumed $\tl{\ph^+}=\ph^-$.
This proves (\ref{T+H}) and the proof of (\ref{W+H}) is
analogous.\sp

Suppose symbols $\ph(z)$ and $\s(\x)$ are related by
\be\s(\x)=\ph\left({1+i\x\ov 1-i\x}\right),\ \ \
\ph(z)=\s\left(i{1-z\ov 1+z}\right). \label{sph}\ee If we use the
fact about the Laguerre polynomials $L_j(x)$ that
\[\int_0^\iy e^{-x}\,e^{-i\x x}\,L_j(2x)\,dx=
\left({1-i\x\ov 1+i\x}\right)^j{1\ov 1+i\x},\]
we see that the $j,k$ entry of the matrix for $W(\s)$ with respect to the orthonormal basis
$f_j(x)=\sqrt2\,e^{-x}\,L_j(2x)$ for $L^2(0,\iy)$ is equal to
\[{1\ov 2\pi}\int_{-\iy}^{\iy}\int_0^{\iy}\int_0^{\iy}
e^{-i\x(x-y)}\,\s(\x)\,f_k(y)\,f_j(x)\,dy\,dx\,d\x\]
\[={1\ov 2\pi}\int_{-\iy}^{\iy} \ph(z(\x))\,z(\x)^{k-j}\,{2\,d\x\ov 1+\x^2}=
{1\ov 2\pi}\int_{|z|=1}\ph(z)\,z^{k-j}\,|dz|.\]
(In the integral on the left $z(\x)$ denotes $(1+i\x)/(1-i\x)$.) Thus the matrix for
$W(\s)$ with respect to this basis is $T(\ph)$. (This fact is implicitly contained
in \cite{R}. Also, of course, the matrix for $W(\s\inv)$ with respect to this basis is $T(\ph\inv)$.)
Similarly the $j,k$ entry of the matrix for $H(\s)$ with
respect to the same basis is equal to
\[{1\ov 2\pi}\int_{-\iy}^{\iy}\int_0^{\iy}\int_0^{\iy}
e^{-i\x(x+y)}\,\s(\x)\,f_k(y)\,f_j(x)\,dy\,dx\,d\x\]
\[={1\ov 2\pi}\int_{-\iy}^{\iy} \ph(z(\x))\,z(\x)^{-j-k}\,{2\,d\x\ov (1+i\x)^2}
={1\ov 2\pi}\int_{-\iy}^{\iy} \ph(z(\x))\,z(\x)^{-j-k-2}\,{2\,d\x\ov (1-i\x)^2}\]
\[={1\ov 2\pi}\int_{|z|=1}\ph(z)\,z^{-j-k-1}\,|dz|.\]
Thus the matrix for $H(\s)$ with respect to this basis is
$H(\ph)$. It follows that the first determinants on the right
sides of (\ref{T+H}) and (\ref{W+H}) are equal and we have
\be{\det\,(W_R(\s)+H_R(\s))\ov\det\,(T_n(\ph)+H_n(\ph))}={G[\s]^R\ov
G[\ph]^n}\; {\det\,(I+Q_R\,H(\s_-/\s_+)\,Q_R)\ov
\det\,(I+Q_n\,H(\ph_-/\ph_+)\,Q_n)}. \label{detratio}\ee In our
application (\ref{sph}) will hold up to constant factors. This is
enough for (\ref{detratio}) to hold since those first determinants
are unchanged if we multiply the symbols by constants.

\begin{center}{\bf III. Fredholm determinant representations in the two-jump case}\end{center}

Our Wiener-Hopf symbol is
\[\s({\x})=\left\{\begin{array}{ll} 1&{\rm if}\ |\x|>1,\\&\\e^{-2\pi i\al}&{\rm if}\ |\x|<1.
\end{array}\right.\]
This can be written alternatively as
\[\s(\x)=\left({\x-1-0i\ov \x-1+0i}\ {\x+1+0i\ov \x+1-0i}\right)^{\al},\]
where the arguments of $\x\pm1\pm i0$ are zero when $\x>0$. The obvious way to
regularize this is to replace 0 by $\ep>0$ everywhere and then
eventually to take the $\ep\to0$ limit. We use a different
regularization which is more complicated at first but will be simpler in the end
after we use the linear-fractional transformation
\be\x(x)=i\,{1-ix\ov 1+ix},\label{xxi}\ee
which maps the unit interval $[-1,\,1]$ to the upper half of the unit circle in the
complex $\x$-plane. To regularize we introduce a parameter $r\in (0,\,1)$ and define
\[\s_r(\x)=\left({\x-\x(r)\ov \x-\overline{\x(r)}}\ {\x+\x(r)\ov \x+
\overline{\x(r)}}\right)^{\al},\]
where the arguments of the four factors are close to zero when $\x$ is large and positive.
Since
\be\x(r)=1+i(1-r)+O((1-r)^2)\ \ \ {\rm as}\ \ r\to 1\label{xiasym}\ee
we see that $\s_r(\x)\to\s(\x)$ as $r\to 1$ when $\x\neq\pm1$.

The Wiener-Hopf factors of $\s_r(\x)$ are
\[\s_r^+(\x)=\left({\x+\x(r)\ov \x-\overline{\x(r)}}\right)^{\al},\ \ \
\s_r^-(\x)=\left({\x-\x(r)\ov \x+\overline{\x(r)}}\right)^{\al},\]
so that
\[{\s_r^+(\x)\ov \s_r^-(\x)}=\left({(\x+\x(r))\,(\x+\overline{\x(r)})
\ov (\x-\x(r))\,(\x-\overline{\x(r)})}\right)^\al.\]
We shall think of $Q_R\,H(\s^-/\s^+)\,Q_R$ as acting on $L^2(0,\,\iy)$ rather than
$L^2(R,\,\iy)$ by making the variable changes $s\to R+s,\ t\to R+t$. If we
observe that $\s_r^\mp(\x)=\s_r^\pm(-\x)$ we see the operator has
kernel
\be K_{R,r}(s,t)={1\ov2\pi}\int_{-\iy}^\iy\left[{\s_r^+(\x)\ov \s_r^-(\x)}-1\right]\,
e^{i\x(2R+s+t)}\,d\x.\label{Hankelkernel}\ee
We deform the contour to the arc of the unit circle in the upper half-plane between
$\x(r)$ and $-\overline{\x(r)}=\x(-r)$ described in two directions (the ``lower'' part
to the right and the ``upper'' part to the left). Then we are going to make the change of
variable (\ref{xxi}) so that the arc of the unit circle corresponds to the $x$-interval
$(-r,\,r)$. A computation shows that
\[{(\x(x)+\x(r))\,(\x(x)+\overline{\x(r)})\ov (\x(x)-\x(r))\,(\x(x)-\overline{\x(r)})}
=-{(1+rx)\,(r+x)\ov(1-rx)\,(r-x)},\]
so that in particular the left side is purely negative on the arc. To compute its argument
on the two sides of the arc, which will be independent of $r$ by continuity, we can
take $r=1$ where the geometry is simple, and find that the argument equals $-\pi$ on
the upper part of the arc and, therefore, $\pi$ on the lower part. Thus
\[{\s_r^+(\x(x))\ov \s_r^-(\x(x))}=e^{-i\pi\al}\,\left({(1+rx)\,(r+x)\ov(1-rx)\,(r-x)}
\right)^\al\]
on the upper part of the arc and
\[{\s_r^+(\x(x))\ov \s_r^-(\x(x))}=e^{i\pi\al}\,\left({(1+rx)\,(r+x)\ov(1-rx)\,(r-x)}
\right)^\al\]
on the lower. Using these we find that (\ref{Hankelkernel}) is equal to
\[i{\sin\pi\al\ov\pi}\int_{-r}^r\left({(1+rx)\,(r+x)\ov(1-rx)\,(r-x)}\right)^\al\,
e^{i\x(x)\,(2R+s+t)}\,{2\ov(1+ix)^2}\,dx.\]

The operator $K_{R,r}$ with this kernel equals a product $UV$ where
$U:L^2(-r,\,r)\to L^2(0,\,\iy)$ has kernel
\[U(s,x)=i{\sin\pi\al\ov\pi}\,e^{i\x(x)s}\,{1\ov 1+ix}\]
and $V:L^2(0,\iy)\to L^2(-r,\,r)$ has kernel
\[V(x,s)=\left({(r+x)\,(1+rx)\ov(r-x)\,(1-rx)}\right)^\al\,{2\ov 1+ix}\,e^{i\x(x)\,(2R+s)}.\]
The determinant of $I+K_{R,r}$ is unchanged if $K_{R,r}$ is replaced by $VU:L^2(-r,\,r)\to
L^2(-r,\,r)$, which has kernel
\[i{\sin\pi\al\ov\pi}\,e^{2i\x(x)R}\,\left({(1+rx)\,(r+x)\ov(1-rx)\,(r-x)}\right)^\al\,
{2\ov(1+ix)\,(1+iy)}\,{i\ov \x(x)+\x(y)}\]
\[=i{\sin\pi\al\ov\pi}\left({(1+rx)\,(r+x)\ov(1-rx)\,(r-x)}\right)^\al\,
{e^{2i\x(x)R}\ov 1+xy}.\]
The determinants will also be unchanged if we change the kernel to
\be i{\sin\pi\al\ov\pi}\left({(1+rx)\,(r+x)\ov(1-rx)\,(r-x)}\,{(1+ry)\,(r+y)\ov(1-ry)\,(r-y)}
\right)^{\al/2}\,{e^{i\,(\x(x)+\x(y))\,R}\ov 1+xy}.\label{WHkernel}\ee\sp

Next, before we take the $r\to1$ limit, we go to Toeplitz and the symbol
\[\ph(e^{i\th})=\left\{\begin{array}{ll}e^{-i\pi\al}&{\rm if}\ -{\pi\ov2}<\th<{\pi\ov2},\\
&\\e^{i\pi\al}&{\rm if}\ {\pi\ov2}<\th<{3\pi\ov2}.\end{array}\right.\]
Thus $\ph(z)=e^{-i\pi\al}$ on the right half of the unit circle and
$\ph(z)=e^{i\pi\al}$ on the left half. This can also be written as
\[\ph(z)=\left({1-iz\ov1+iz\inv}\,{1-iz\inv\ov 1+iz}\right)^\al,\]
where the arguments are chosen so that the factors involving $z$ have argument 0 at $z=0$ and
those involving $z\inv$ have argument 0 at $z=\iy$. To regularize we choose the same
parameter $r$ as before and replace $\ph(z)$ by
\[\ph_r(z)=\left({1-irz\ov 1+irz\inv}\,{1-irz\inv\ov 1+irz}\right)^\al.\]

We point out here that if $z(\x)=(1+i\x)/(1-i\x)$ then a computation shows that
\[\ph_r(z(\x))=\left({1+ir\ov 1-ir}\right)^{2\al}\s_r(\x).\]
Thus (\ref{sph}) holds for the symbols $\s_r(\x)$ and $\ph_r(z)$, up to constant
factors, and therefore we have the identity (\ref{detratio}).

To continue, the Wiener-Hopf factors of $\ph_r$ are
\[\ph_r^+(z)=\left({1-irz\ov 1+irz}\right)^\al,\ \ \
\ph_r^-(z)=\left({1-irz\inv\ov 1+irz\inv}\right)^\al,\]
so that
\[{\ph_r^+(z)\ov \ph_r^-(z)}=\left({(1-irz)\,(1+irz\inv)\ov(1+irz)\,(1-irz\inv)}\right)^\al.\]

If we think of $Q_n\,H(\ph^-/\ph^+)\,Q_n$ as acting on $\ell^2(\Zp)$
then it has $j,k$ entry
\[H_{n,r}(j,k)={1\ov2\pi i}\int{\ph_r^+(z)\ov \ph_r^-(z)}\,z^{2n+j+k}\,dz\]
with integration over the unit circle. (We use here the fact that $\ph_r^-(z)=\ph_r^+(z\inv)$.)
With the variable change $z\to iz$ the above becomes
\[{i^{2n+j+k}\ov2\pi}\int\left({(1+rz)\,(1+rz\inv)\ov(1-rz)\,(1-rz\inv)}\right)^\al
\,z^{2n+j+k}\,dz.\]

The integrand is analytic inside the
unit circle cut on the line segment $[-r,\,r]$. We deform the path of
integration to this segment described back and forth. The expression in large parentheses
is real and negative on the segment. The limit of its argument from above equals $-\pi$
and from below equals $\pi$. Hence
\[H_{n,r}(j,k)=i^{2n+j+k+1}{\sin\pi\al\ov\pi}\int_{-r}^r\left({(1+rx)\,(r+x)\ov(1-rx)\,(r-x)}
\right)^\al\,x^{2n+j+k}\,dx.\]
This operator $H_{n,r}$ is now a product $UV$ where $U:L^2(-r,r)\to\ell^2(\Zp)$
has kernel
\[U(j,x)=i^{2n+j+1}{\sin\pi\al\ov\pi}\,x^j\]
and $V:\ell^2(\Zp)\to L^2(-r,r)$ has kernel
\[V(x,j)=i^j\,\left({(1+rx)\,(r+x)\ov(1-rx)\,(r-x)}\right)^\al\,x^{2n+j}.\]
For $\det(I+H_{n,r})$ this may be replaced by $VU:L^2(-r,r)\to L^2(-r,r)$ which has
kernel
\[\sum_{j=0}^\iy V(x,j)\,U(j,y)=i{\sin\pi\al\ov\pi}
\left({(1+rx)\,(r+x)\ov(1-rx)\,(r-x)}\right)^\al{(ix)^{2n}\ov 1+xy}.\]
Without affecting the determinants this can be changed to
\be i{\sin\pi\al\ov\pi}
\left({(1+rx)\,(r+x)\ov(1-rx)\,(r-x)}\,{(1+ry)\,(r+y)\ov(1-ry)\,(r-y)}\right)^{\al/2}
\,{(-xy)^n\ov 1+xy}.\label{Tkernel}\ee

We change notation and now denote by $K_{R,r}(x,y)$ and
$H_{n,r}(x,y)$ the kernels (\ref{WHkernel}) and
(\ref{Tkernel}), respectively. They are independent of $r$ but the corresponding 
operators both act on $L^2(-r,r)$ and we have shown that
\[{\det\,(W_R(\s_r)+H_R(\s_r))\ov\det\,(T_n(\ph_r)+H_n(\ph_r))}={G[\s_r]^R\ov G[\ph_r]^n}\,
{\det\,(I+K_{R,r})\ov\det\,(I+H_{n,r})}.\]

Suppose we knew that $I+H_{n,r}$ was invertible. Then we could 
rewrite the above as
\[{\det\,(W_R(\s_r)+H_R(\s_r))\ov\det\,(T_n(\ph_r)+H_n(\ph_r))}=
{G[\s_r]^R\ov G[\ph_r]^n}\,
\det\,\left((I+K_{R,r})\,(I+H_{n,r})\inv\right).\]
Denote by $K_R^0$ and $H_n^0$ the corresponding operators on $L^2(-1,1)$. We shall
show in the next section that $I+H_n^0$ is uniformly invertible for large
$n$ and a trivial modification shows that $I+H_{n,r}$ is
uniformly invertible if in addition $r$ is bounded away from zero. We shall also show that
$K_R^0-H_n^0$ is trace class. Granting these things in advance, the above holds and in addition
we deduce the limiting relation
\be{\det\,(W_R(\s)+H_R(\s))\ov\det\,(T_n(\ph)+H_n(\ph))}=e^{-2i\al R}\;
\det\,\left((I+K_R^0)\,(I+H_n^0)\inv\right).\label{detquotient}\ee
We used here that $G[\ph]=1$ and $G[\s]=e^{-2i\al}$.
\sp

\begin{center}{\bf IV. Asymptotics}\end{center}\sp

Theorem 2.4 of \cite{BE2} states for our $T_n(\ph)+H_n(\ph)$  that
\[\det\,(T_n(\ph)+H_n(\ph))\sim n^{-3\,\al^2}\,2^{4\,\al^2}\,G(1-2\al)\,G(1+2\al)\]
as $n\to\iy$. Thus our Theorem will be proved if we can show that
(\ref{detquotient}) is valid and that the determinant on the right side there
tends to 1 as $n$ tends to infinity in such a way that $n\sim R$. These will
follow once we show that  $K_R^0-H_n^0$ tends to 0 in trace norm and that the
operators $I+H_n^0$ are uniformly invertible.\footnote
{In fact we shall show that the latter
holds for all $\al$ satisfying $|\Re\,\al|<1/2$ except for those lying in a 
discrete set. The extra condition on $\al$ can be removed at the end by an easy 
analyticity argument as at the end of \cite{BW}. We shall also assume that
$R-n=o(n^{1/2})$. This is good enough since any $n=n(R)$ satisfying $n(R)\sim R$ as $R\to\iy$
would give us the asymptotics.}

For the first, we observe that the ``problem points'' for our kernels 
are the pairs \linebreak $x=\pm1,\ y=\mp1$. More precisely, for any $\dl>0$ the operator with kernel
\[{1-[\ch_{(1-\dl,\,1)}(x)\,\ch_{(-1,\,-1+\dl)}(y)+
\ch_{(-1,\,-1+\dl)}(x)\,\ch_{(1-\dl,\,1)}(y)]\ov 1+xy}\]
is trace class, and the same is true even with the factor
\[\left({1+x\ov 1-x}\,{1+y\ov1-y}\right)^{\al}\]
long as $|\Re\,\al|<1/2$. Moreover multiplication by $x^n$ converges strongly to 0
as $n\to\iy$. It follows that
\[\Big(1-[\ch_{(1-\dl,\,1)}(x)\,\ch_{(-1,\,-1+\dl)}(y)+
\ch_{(-1,\,-1+\dl)}(x)\,\ch_{(1-\dl,\,1)}(y)]\Big)\,H_n^0(x,y)=o_1(1).\]
(By this we mean that the left side is the kernel of an operator
whose trace norm is $o(1)$.) A similar argument applies to
$K_R^0$, so we are left with showing that the kernels
\[\ch_{(1-\dl,\,1)}(x)\,\ch_{(-1,\,-1+\dl)}(y)\,(K_R^0(x,y)-H_n^0(x,y)),\ \ \
\ch_{(-1,\,-1+\dl)}(x)\,\ch_{(1-\dl,\,1)}(y)\,(K_R^0(x,y)-H_n^0(x,y))\]
are $o_1(1)$. We consider only the first, and it is convenient to make the variable
change $y\to-y$ so that it becomes a constant times
\be\ch_{(1-\dl,\,1)}(x)\,\ch_{(1-\dl,\,1)}(y)\,\left({1+x\ov 1-x}\,{1-y\ov1+y}\right)^{\al}
{e^{i\,(\x(x)+\x(-y))\,R}-(xy)^n\ov 1-xy}.\label{K-H}\ee

Lemma 1 of \cite{BW} says, with slightly different notation, that the trace norm of a kernel
$f(\x)g(\e)/(\x+\e)$ on
$L^2(0,\,\iy)$ is at most a constant depending on $b$ times the square root of
\[\int|f(\x)|^2{d\x\ov \x^{1+b}}\,\cdot\,\int|g(\e)|^2{d\e\ov \e^{1-b}}.\]
Here $b$ belongs to $(-1,1)$ but is otherwise arbitrary. If we make the
substitutions
\[\x={1-x\ov 1+x},\ \ \ \eta={1-y\ov1+y}\]
we find that the trace norm of a kernel $F(x)G(y)/(1-xy)$ on $L^2(0,1)$ is at
most a constant times the square root of
\be\int_0^1(1-x)^{-b-1}\,(1+x)^{b-1}\,|F(x)|^2\,dx\,\cdot\,
\int_0^1(1-y)^{b-1}\,(1+y)^{-b-1}\,|G(y)|^2\,dy.\label{intprod}\ee

If we write
\[e^{i\,(\x(x)+\x(-y))\,R}-(xy)^n\]
\be=\left(e^{i\,(\x(x)-1)R}-x^n\right)\,
e^{i\,(\x(-y)+1)\,R}+x^n\left(e^{i\,(\x(-y)+1)\,R}-y^n\right),\label{kernsum}\ee
then (\ref{K-H}) is written correspondingly as the sum of two operators of the above
form. In the first
\[F(x)=\ch_{(1-\dl,\,1)}(x)\,\left({1+x\ov 1-x}\right)^{\al}\,
\left(e^{i\,(\x(x)-1)R}-x^n\right),\]
\[G(y)=\ch_{(1-\dl,\,1)}(y)\,\left({1-y\ov 1+y}\right)^{\al}\,e^{i\,(\x(-y)+1)\,R}.\]

The second integral in (\ref{intprod}) is $O(1)$ as long as $b>-2\,\Re\,\al$. (Recall
that $\Im\,\x(-y)>0$.)  We choose $b=-2\,\Re\,\al+\ep$ with $\ep>0$ small enough so that
the condition $b<1$ is satisfied. We are left with  showing that
\[\int_{1-\dl}^1\,(1-x)^{-1-\ep}\,|e^{i\,(\x(x)-1)R}-x^n|^2\,dx=o(1).\]
It is enough to show this with $\ep$ replaced by $1/2$, say. Using (\ref{xiasym}) and the
variable change $x\to 1-x$ we see that we want to show
\[\int_0^\dl x^{-3/2}\,|e^{(-x+O(x^2))\,R}-e^{(-x+O(x^2))\,n}|^2\,dx=o(1).\]
The integral over $x>n^{-2/3}$ is exponentially small (if $\dl$ was chosen small
enough so that the $O(x^2)$ terms are less than $x/2$) and the rest is at most
\[\int_0^{n^{-2/3}}x^{-3/2}\,|x\,(n-R)+O(x^2n)|^2\,dx=O(n\inv(n-R)^2+n^{-1/3}),\]
and this is $o(1)$ if $R-n=o(n^{1/2})$.

For the operator arising from the second summand in (\ref{kernsum}), in particular its
second factor, we
use the fact that $\x(-y)=-\overline{\x(y)}$ and then (\ref{xiasym}) again. We choose
$b=\Re\,\al+1/2$ in (\ref{intprod}) and the same
integrals arise as before. This completes the proof that $K_R^0-H_n^0=o_1(1)$.

It remains to establish the uniform invertibility of $I+H_n^0$ for all $\al$ satisfying
$|\Re\,\al|<1/2$ except for those lying in a discrete set. Recall that we only change the
operator by $o_1(1)$ if we multiply the kernel by
\[\ch_{(0,\,1)}(x)\,\ch_{(-1,\,0)}(y)+\ch_{(-1,\,0)}(x)\,\ch_{(0,\,1)}(y)\]
so this will not affect uniform invertibility. Think of $H_n^0$ as already having this
factor and then think of $L^2(-1,1)$ as
$L^2(0,1)\oplus L^2(-1,0)$. The new $H_n^0$ will have a certain matrix representation.
If we make
the variable change $x\to-x$ in $L^2(-1,0)$ so that it becomes $L^2(0,1)$, then all
operators in the matrix act on $L^2(0,1)$, and the matrix kernel for $H_n^0$ becomes
\[i{\sin\pi\al\ov\pi}\twotwo{0}{{\hspace{-4ex}}\left({1+x\ov 1-x}\,{1-y\ov1+y}
\right)^{\al}\,{(xy)^n\ov 1+xy}}{\left({1-x\ov 1+x}\,{1+y\ov1-y}
\right)^{\al}\,{(xy)^n\ov 1+xy}}{{\hspace{-4ex}}0}.\]

First we make the same variable changes
\[x\to {1-x\ov 1+x},\ \ \ y\to {1-y\ov 1+y}\]
as before, and the matrix kernel becomes
\[i{\sin\pi\al\ov\pi}\twotwo{0}{{\hspace{-4ex}}\left({1-x\ov 1+x}\,{1-y\ov1+y}
\right)^n{(y/x)^{\al}\ov x+y}}{\left({1-x\ov 1+x}\,{1-y\ov1+y}
\right)^n{(x/y)^{\al}\ov x+y}}{{\hspace{-4ex}}0}\]
still acting on $L^2(0,1)\oplus L^2(0,1)$. Then, using the lemma quoted from \cite{BW}
and computations analogous to those already done,
we can show that the error incurred is $o_1(1)$ if we replace
\[\left({1-x\ov 1+x}\,{1-y\ov1+y}\right)^n\]
by $e^{-2n(x+y)}$, so this will not affect uniform invertibility. Then we make the
substitutions $x\to x/2n,\ \ y\to y/2n$ and the kernel becomes
\[i{\sin\pi\al\ov\pi}\twotwo{0}{{\hspace{-2ex}}e^{-(x+y)}\,{(y/x)^{\al}\ov x+y}}
{e^{-(x+y)}\,{(x/y)^{\al}\ov x+y}}{{\hspace{-2ex}}0}\]
on $L^2(0,2n)\oplus L^2(0,2n)$.

This is a variant of the operator that appears in Lemma~3.2 
of \cite{BW}. The argument is essentially the same here, and we
outline it. It suffices to show that $I$ plus the limiting operator on $L^2(0,\iy)\oplus L^2(0,\iy)$
is invertible, except
for a discrete set of $\al$. Call the kernel $L(x,y)$. First we show, by transforming to a
Wiener-Hopf operator with matrix kernel, that $I$ plus the operator with kernel
\[L_0(x,y)=i{\sin\pi\al\ov\pi}\,\ch_{(0,1)}(x)\,\ch_{(0,1)}(y)\,
\twotwo{0}{{(y/x)^\al\ov x+y}}{(x/y)^\al\ov x+y}{0}\]
is invertible except for a discrete set of $\al$. Then we observe that
$L-L_0$ is a trace class operator. Putting these together we see
that the invertibility of $I+L$ is equivalent to the
invertibility of $I+(I+L_0)\inv (L-L_0)$, which in turn is equivalent to the
nonvanishing of its
determinant. We know this does not vanish if $\al$ is sufficiently small and so
it can only vanish for a discrete set of $\al$.\sp

\begin{center}{\bf Acknowledgments}\end{center}

The first author was supported by National Science Foundation grant DMS-9970879
and the third author by grant DMS-9732687.

\end{document}